\documentclass[11pt]{amsart}
\usepackage{latexsym}
\usepackage{amssymb}
\usepackage{amscd}
\usepackage[mathscr]{euscript}

\normalsize

\RequirePackage{xspace}

\newcommand{\thought}[1]{}
\renewcommand{\thought}[1]{ \textbf{[#1]}}

\usepackage{enumerate}        

\newcommand\nc {\newcommand}
\newcommand\rnc{\renewcommand}
\nc\script{\mathscr}

\newcount\blopone
\newcount\xone
\newcount\xtwo
\newcount\ytwo

\newtheorem{theorem}{Theorem}[section]
\newtheorem{prop}[theorem]{Proposition}

\newtheorem{refinement}[theorem]{Refinement}
\newtheorem{summary}[theorem]{Summary}
\newtheorem{importnota}[theorem]{Important Notation}
\newtheorem{prblm}[theorem]{Problem}
\newtheorem{notation}[theorem]{Notation}
\newtheorem{defin}[theorem]{Definition}
\newtheorem{caution}[theorem]{Caution}
\newtheorem{remark}[theorem]{Remark}
\newtheorem{reminder}[theorem]{Reminder}
\newtheorem{lemma}[theorem]{Lemma}
\newtheorem{construction}[theorem]{Construction}
\newtheorem{corollary}[theorem]{Corollary}
\newtheorem{example}[theorem]{Example}
\newtheorem{conclusion}[theorem]{Conclusion}
\newtheorem{triviality}[theorem]{Triviality}
\newtheorem{proto}[theorem]{Prototype Quasifibration}
\newtheorem{cauex}[theorem]{Cautionary Example}
\newtheorem{hypo}[theorem]{Hypothesis}
\newtheorem{subth}{ }[theorem]
\newtheorem{case}{Case}[theorem]
\newtheorem{ssubth}{ }[subth]

\nc\tri[1]{\begin{triviality}
\label{#1}}
\nc\cas[1]{\begin{case}
\label{#1}
\begin{em}}
\nc\rfn[1]{\begin{refinement}
\label{#1}}
\nc\prt[1]{\begin{proto}
\label{#1}}
\nc\lem[1]{\begin{lemma}
\label{#1}}
\nc\pro[1]{\begin{prop}
\label{#1}}
\nc\thm[1]{\begin{theorem}
\label{#1}}
\nc\cor[1]{\begin{corollary}
\label{#1}}
\nc\dfn[1]{\begin{defin}
\begin{em}
\label{#1}}
\nc\sthm[1]{\begin{subth}
\label{#1}}
\nc\exm[1]{\begin{example}
\label{#1}
\begin{em}}
\nc\plm[1]{\begin{prblm}
\label{#1}
\begin{em}}
\nc\rmk[1]{\begin{remark}
\label{#1}
\begin{em}}
\nc\rmd[1]{\begin{reminder}
\label{#1}
\begin{em}}
\nc\ntn[1]{\begin{notation}
\label{#1}
\begin{em}}
\nc\smr[1]{\begin{summary}
\label{#1}
\begin{em}}
\nc\cau[1]{\begin{caution}
\label{#1}
\begin{em}}
\nc\hyp[1]{\begin{hypo}
\label{#1}}
\nc\imn[1]{\begin{importnota}
\label{#1}
\begin{em}}
\nc\cax[1]{\begin{cauex}
\label{#1}
\begin{em}}
\nc\con[1]{\begin{construction}
\label{#1}
\begin{em}}
\nc\ssthm[1]{\begin{ssubth}
\label{#1}
\begin{em}}
\nc\cnc[1]{\begin{conclusion}
\label{#1}
\begin{em}}

\nc\elem{\end{lemma}}
\nc\erfn{\end{refinement}}
\nc\eprt{\end{proto}}
\nc\ethm{\end{theorem}}
\nc\ecor{\end{corollary}}
\nc\edfn{$\,$\hfill{$\Box$}\vskip 2mm\end{em}
\end{defin}}
\nc\esthm{\end{subth}}
\nc\epro{\end{prop}}
\nc\etri{\end{triviality}}
\nc\eexm{$\,$\hfill{$\Box$}\vskip 2mm\end{em}
\end{example}}
\nc\ermk{$\,$\hfill{$\Box$}\vskip 2mm\end{em}
\end{remark}}
\nc\ermd{$\,$\hfill{$\Box$}\vskip 2mm\end{em}
\end{reminder}}
\nc\eplm{\end{em}
\end{prblm}}
\nc\ecas{\end{em}
\end{case}}
\nc\ecau{\end{em}
\end{caution}}
\nc\ecax{\end{em}
\end{cauex}}
\nc\eimn{\end{em}
\end{importnota}}
\nc\entn{\hfill{$\Box$}\end{em}
\end{notation}}
\nc\econ{\end{em}
\end{construction}}
\nc\esmr{\end{em}
\end{summary}}
\nc\ehyp{
\end{hypo}}
\nc\ecnc{\end{em}
\end{conclusion}}
\nc\essthm{\end{em}
\end{ssubth}}

\newtheorem{proposition}[theorem]{Proposition}

\nc\sst{\scriptstyle}
\newcommand{\comment}[1]{}
\newcommand{\ri}{\longrightarrow}

\nc\bR{{\mathbf R}}
\nc\bS{{\mathbf S}}
\nc\bT{{\mathbf T}}
\nc\bU{{\mathbf U}}
\nc\z{\zeta}
\nc\bc{{\mathbb{BC}}}
\nc\ck{{\script K}}
\nc\ct{{\script T}}
\nc\cs{{\script S}}
\nc\car{{\script R}}
\nc\ca{{\script A}}
\nc\cb{{\script B}}
\nc\cc{{\script C}}
\nc\cd{{\script D}}
\nc\ce{{\script E}}
\nc\ci{{\script I}}
\nc\co{{\script O}}
\nc\bZ{{\mathbb Z}}
\nc\bd{\begin{description}}
\nc\ed{\end{description}}
\nc\ctob{{\script C}at\big(\ci^{op},\ca\big)}
\nc\clim{{\ds\mathop{\rm lim}_{\ds\longleftarrow}}}
\nc\climi{\clim^{\!i}\,}
\nc\climn{\clim^{\!n}\,}
\nc\colim{{\ds\mathop{\rm colim}_{\ds\la}}}
\nc\oa{\overline{\ca}}
\nc\s{\sigma}
\nc\ta{\tau}
\nc\os{\overline\sigma}
\nc\ot{\overline\tau}
\nc\T{\Sigma}
\nc\de[1]{{\mathop{\rm deg(#1)}}}
\nc\Ad[1]{\mathop{\rm Ad}(#1)}
\nc\ad[1]{\mathop{\rm ad}(#1)}

\def\der #1 {D\left(#1\right)}
\nc\prf{\begin{proof}}
\nc\eprf{\end{proof}}
\nc\ds{\displaystyle}

\nc\ab{{\script A}b}

\nc\csab{{\script C}at\big(\cs^{op},\ab\big)}
\nc\ctab{{\script C}at\Big({\{\ct^\alpha\}}^{op},\ab\Big)}
\nc\csex{{\script E}x\big(\cs^{op},\ab\big)}
\nc\ctex{{\script E}x\Big({\{\ct^\alpha\}}^{op},\ab\Big)}
\nc\sub{\qquad\subset\qquad}
\nc\ctr[1]{{\left.\ct\left(-,#1\right)\right|}_{\cs}}
\nc\ctrf[2]{{\left.\ct\left(#1,#2\right)\right|}_{\cs}}
\nc\Ctr[1]{{\left.\ct\left(-,#1\right)\right|}_{\ct^\alpha}}
\nc\Ctrf[2]{{\left.\ct\left(#1,#2\right)\right|}_{\ct^\alpha}}

\nc\la{\longrightarrow}
\nc\oti{{^L\otimes_R^{}}}
\nc\rs{\s^{-1}R}
\nc\br{{\{\s^{-1}R\}}}

\nc\nin{\noindent}
\nc\cad[1]{\text{card}(#1)}
\nc\eq{\quad=\quad}
\nc\BA{\begin{array}{c}}
\nc\EA{\end{array}}
\nc\kth{{\it K}--theory}

\nc\barr{
\[
\begin{array}{cccccccccccccccc}
}
\nc\earr{
\end{array}
\]
}
\nc\as[1]{{\langle S\rangle}^{#1}}
\nc\sh{\hbox{\it shift}}

\nc\yy[1]{{\left.\ct\left(-,#1\right)\right|}_{\ct^c}}
\nc\vrep[2]{{\left.\ct\left(#1,#2\right)\right|}_{\ct^\alpha}}
\nc\da{\downarrow}
\nc\Hom{{\mathop{\rm Hom}}}
\nc\End{{\mathop{\rm End}}}
\nc\Ext{{\mathop{\rm Ext}}}
\nc\Tor{{\mathop{\rm Tor}}}
\nc\mr\Modtc
\nc\PExt{{\mathop{\rm PExt}}}
\nc\bA{{\mathbf A}}
\nc\bB{{\mathbf B}}
\nc\bC{{\mathbf C}}
\nc\bD{{\mathbf D}}

\nc\ra\ri
\nc\y[1]{\mathbf{y}#1}
\nc\x[1]{\mathbf{z}#1}
\nc\Mod[1]{\ensuremath{\mathop{\textup{Mod-}#1}}\xspace}
\nc\Md {\ensuremath{\mathop{\textup{Mod}}}}
\rnc\mod[1]{\ensuremath{\mathop{\textup{mod-}#1}}\xspace}

\nc\Modtc{\Mod{\ct^c}}
\nc\pgldim[1]{\mathop{\rm pgldim}\,#1}

\begin{document}

\author{Amnon Neeman}
\address{A.N. : Center for Mathematics and its Applications \newline
\indent School of Mathematical Sciences\newline
\indent John Dedman Building\newline
\indent The Australian National University\newline
\indent Canberra, ACT 0200\newline
\indent AUSTRALIA}
\email{Amnon.Neeman@anu.edu.au}
\author{Andrew Ranicki}
\address{A.R. : Department of Mathematics and Statistics \newline
\indent University of Edinburgh \newline
\indent Edinburgh EH9 3JZ\newline
\indent SCOTLAND, UK}
\email{aar@maths.ed.ac.uk}
\author{Aidan Schofield}
\address{A.S. : School of Mathematics\newline
\indent University of Bristol\newline
\indent Bristol BS8 1TW\newline
\indent ENGLAND, UK}
\email{Aidan.Schofield@bristol.ac.uk}

\title[Representations of algebras as universal localizations]
{Representations of algebras as universal localizations}

\date{July 2, 2002}


\keywords{representations of algebras, noncommutative localization\\
\indent A.R. acknowledges the support of the Leverhulme Trust}
\maketitle

\tableofcontents

\section*{Introduction}
\label{S0}

Given a presentation of a finitely presented group, there is a natural
way to represent the group as the fundamental group of a
$2$-complex. The first part of this paper demonstrates one possible
way to represent a finitely presented algebra $S$ in a similarly
compact form.  From a presentation of the algebra, we construct a
quiver with relations whose path algebra is finite dimensional. When
we adjoin inverses to some of the arrows in the quiver, we show that
the path algebra of the new quiver with relations is $M_n(S)$ where
$n$ is the number of vertices in our quiver. The slogan would be that
every finitely presented algebra is Morita equivalent to a universal
localization of a finite dimensional algebra.

Two applications of this are then considered. Firstly, given a ring
homomorphism $A \to B$, we say that $B$ is {\it stably
flat} over $A$ if and only if $\Tor_i^A(B,B)=0$ for all $i>0$.  In a
recent paper \cite{Neeman-Ranicki01}, the first two authors show that
there is a long exact sequence in algebraic $K$-theory associated to a
universal localization provided the localization is stably flat.  
Given a finitely presented algebra $S$ we construct in Section 1 a
finite dimensional algebra $A$ with a universal localization
$\sigma^{-1}A=M_n(S)$ (Theorem \ref{t0}), such that $A$ has global dimension 
$\leqslant 2$ (Proposition \ref{p0}).  For a stably flat universal localization
$A \to \sigma^{-1}A$ the global dimensions are such that
${\rm g.d.}(\sigma^{-1}A) \leqslant {\rm g.d.}(A)$ (Lemma \ref{l0}). 
Thus any finitely presented algebra $S$ of global dimension $\geqslant 3$ 
provides an example of a universal localization which is not stably flat.

Secondly, the Malcolmson normal form states that every element of the
localised ring can be written in the form $as^{-1}b$ where $s\colon
P\rightarrow Q$ lies in the upper triangular closure of $\sigma$,
$a\colon A\rightarrow Q$ and $b\colon P\rightarrow A$ are maps in the
category of finitely generated projective modules over the original
ring $A$, and gives an equivalence condition on such elements which
determines when they define the same element of the localised
ring. This equivalence condition depends on the existence of certain
maps in $\sigma$ and the category of finitely generated projective
modules. One might reasonably ask if such an equation could be
constructed algorithmically. We show that this cannot be done.

In Section 2, we consider a related construction of a ring by
universal localization where we calculate explicitly the values of
$\Tor_i^A({\sigma}^{-1}A,{\sigma}^{-1}A)$. For any $n \geqslant 3$
we obtain an injective universal localization $A \to \sigma^{-1}A$ with
$\Tor_i^A({\sigma}^{-1}A,{\sigma}^{-1}A)=0$ for $1 \leqslant i \leqslant n-2$
and $\neq 0$ for $i=n-1$.

\section{Algebras}
\label{S1}

An {\it algebra} over a field $k$ is a ring $A$ with a homomorphism from $k$ to the centre of $R$.
By definition, the algebra $A$ is {\it finite dimensional} if it is
a finite dimensional vector space over $k$.
By definition, an algebra $S$ is {\it finitely presented} if
it has a finite number of generators and relations, so that it has the form
$$S~=~k \langle x_1,x_2,\dots,x_a \rangle/\langle y_1,y_2,\dots,y_b\rangle~.$$
A finite dimensional algebra $S$ is finitely presented, since for
any basis $e_1,e_2,\dots,e_a$ the coefficients $c_{pqr} \in k$ in
$$e_pe_q~=~\sum\limits_r c_{pqr}e_r \in S~~(1 \leqslant p,q,r \leqslant  a)$$
are such that
$$S~=~k \langle x_1,x_2,\dots,x_a \rangle/
\langle x_px_q - \sum\limits_r c_{pqr}x_r\rangle~.$$
\indent For any finitely presented algebra $S$ over $k$ we shall exhibit the
matrix algebra $M_n(S)$ for some integer $n$ as the universal
localization $\sigma^{-1}A$ of a finite dimensional algebra $A$ over
$k$ inverting a finite set $\sigma$ of maps between finitely generated projective
$A$-modules. We shall construct $A$ as the path algebra of a quiver with
relations, and it will be clear from the construction that $A$ is of
global dimension $2$, but the natural map $A \to \sigma^{-1}A=M_n(S)$
may not be an injection.  Then a variation
of the construction allows us to ensure that $A\to \sigma^{-1}A$ is
injective and $A$ has finite global dimension.  From this it is fairly
clear that for suitable choice of $S$, for example, of infinite global
dimension, the $\Tor_i^{A}(\sigma^{-1}A,\sigma^{-1}A)$ cannot all
vanish. We present examples to show  variations on these techniques.  \medskip

First of all,  recall the language of quivers with relations.
\medskip

A {\it quiver} $Q$ has a finite {\it vertex set} $V_Q=\{v,w,\dots\}$ and
finite {\it arrow set} $A_Q=\{a,b,\dots\}$.  Each arrow $a \in A$ has a {\it
tail} $ta\in V_Q$ and {\it head} $ha\in V_Q$.  A {\it path of length $i$} is a
formal word in the arrows $a_1,\dots,a_i$ such that for $1\leqslant  j < i$,
$ha_j=ta_{j+1}$.  Its tail is $ta_1=v$ and its head is $ha_i=w$ and we
say that it is a path from $v$ to $w$.  For each vertex $v\in V_Q$ we
have a path $f_v\in A_Q$ of length $0$ at $v$ whose head and tail are
both $v$.  For vertices $v$ and $w$, we define $[v,w]$ to be the vector
space with basis the set of paths from $v$ to $w$.
The {\it path algebra} of $Q$ is the vector space
$$\Lambda(Q)~=~\bigoplus\limits_{v,w\in V_Q} [v,w]$$
with the product given by the composition of arrows, which makes it
into an associative algebra with $1=\sum_v f_v$.  Note that this
composition gives an injective linear map from $[u,v]\otimes [v,w]$ to
$[u,w]$.  \medskip

A {\it quiver with relations}  $(Q,R)$ is a quiver $Q$ together with a set of relations
$R=\{r_i\}$ where each $r_i$ is an element of $\cup_{v,w} [v,w]$. In the
examples we shall be discussing $R$ is a finite set. Each element $r$ of
$R$ has a well-defined head and tail which we shall write as $tr$ and
$hr$. For vertices $v,w$, define $R[v,w]$ to be the linear subspace of
$[v,w]$ of the form $\sum_{r\in R}[v,tr]r[hr,w]$. Then $\oplus_{v,w}
R[v,w]$ is an ideal in $\Lambda(Q)$ and the factor algebra
$$\Lambda(Q,R)~=~\Lambda(Q)/\oplus_{v,w}R[v,w]$$
is called the {\it path algebra} of the quiver with relations
$(Q,R)$. We define
$$(v,w)~=~[v,w]/R[v,w]~,$$
so
$$\Lambda(Q,R)~=~\bigoplus\limits_{v,w\in V_Q} (v,w)~.$$

\medskip

We begin with notation. Let
$$S~=~k\langle X: Y\rangle $$
where $X=\{x_i: 1 \leqslant  i \leqslant  a\}$ and $Y = \{y_j: 1\leqslant
j \leqslant  b\}$ is a finite subset of $k\langle X\rangle $.  In turn,
each element of $Y$ can be written in a unique way as a linear
combination of words in the set $X$.  Thus
$$y_j~=~\sum\limits_{\ell=1}^{c_j}\lambda_{j\ell}w_{j\ell}$$
 for suitable elements $\lambda_{j\ell}\in k$ and words
$w_{j\ell}$.  Let $n-1$ be the maximal length of a word $w_{j\ell}$.
\medskip

We consider the quiver $Q$ with vertex and arrow sets
$$V_Q~=~\{1,\dots ,n\}~,~A_Q~=~\{e_1,\dots ,e_{n-1}\}\cup \{1,\dots,n-1\}\times X$$
where $e_m$ is an arrow from $m$ to $m+1$ and
$a_{mi}=(m,x_i)$ is also an arrow from $m$ to $m+1$. Eventually we
are going to invert the arrows $e_m$ and then they and their
inverses will generate a copy of $M_n(k)$. With this in mind and for
convenience of notation we define for $1\leqslant  s < t \leqslant  n$, $e_{s,t} =
e_{s}\dots e_{t-1}$. Thus $e_{s,t}$ is the unique path using the arrows
$e_m$ from the vertex $s$ to the vertex $t$. We also define
$e_{m,m}$ to be the empty path from $m$ to $m$.
\medskip

We construct a set of relations on this quiver. Our first set of
relations is
$$T~=~\{t_{mi}: 1 < m < n, 1 \leqslant  i \leqslant  a\}$$
where
$$t_{mi}~=~a_{1i}e_{2,n} - e_{1,m}a_{mi}e_{m+1,n}~.$$
These, in a sense
which will become clear soon, ensure that $a_{mi}$ for fixed $i$ all
represent the element $x_i$. Now let $w=x_{i_1}\dots x_{i_{u}}$ be a
word of length less than $n$. We define $w'=a_{1,i_1}\dots
a_{u,i_{u}}e_{u+1,n}$, a path in the quiver $Q$ from $1$ to $n$. We
define
$$Y' ~=~ \{y'_j: 1 \leqslant  j\leqslant  b\}$$
where
$$y'_j~=~\sum_{\ell=1}^{c_j} \lambda_{j\ell}w'_{j\ell}~.$$
Our relations on the
quiver are $T\cup Y'$. Its path algebra $A$ is evidently finite
dimensional and it is a simple matter as we shall see to check that $A$
has global dimension $2$.
\medskip

For each vertex $m$, let $P_m$ be the corresponding projective
representation of the quiver $Q$. Given a path $p$ in the quiver $Q$
from $s$ to $t$, there is a corresponding map $\hat{p}\colon
P_{t}\rightarrow P_{s}$. We shall abuse notation by writing $P_m$ for
the corresponding projective module for $A$ and $\hat{p}$ for the
corresponding homomorphism of projective $A$ modules. Let
$\sigma=\{\hat{e}_1,\dots ,\hat{e}_{n-1}\}$.

\begin{theorem}
\label{t0}
For $S$, $A$, $\sigma$ and $n$ defined as above there is an isomorphism of rings
$$\sigma^{-1}A~\cong~M_n(S)~.$$
\end{theorem}
\begin{proof}
We continue to use the notation from the preceding discussion. 
In the special case when $S=k\langle x_1,x_2,\dots,x_a\rangle$ has no relations, it is
well-known that $M_2(S)=\sigma^{-1}A$ is the universal localization
inverting the top arrow in the path algebra $A$ of the
quiver $Q$ with two vertices and $a+1$ arrows. 
From now on, we assume that $S=k\langle X:Y \rangle$ has at least one
relation, i.e. $b \geqslant 1$. 
We enlarge our quiver with relations $Q$ to a quiver with relations $Q'$
by adjoining arrows $f_m$ from the vertex $m+1$ to the vertex $m$
together with relations $e_mf_m=e_{m,m}$ and
$f_me_m=e_{m+1,m+1}$ for $1 \leqslant  m \leqslant  n-1$. The path algebra
$A'$ of the quiver with relations $Q'$ is just $\sigma^{-1}A$. If we
consider the subquiver with relations with the same vertex set and
arrows $e_m,f_m$ for $1 \leqslant  m \leqslant  n-1$, it is clear that there
is a unique path $e_{s,t}$ from vertex $s$ to vertex $t$ that involves
no subpath of type $e_mf_m$ or of type $f_me_m$ and that
$e_{s,t}e_{t,u}=e_{s,u}$ for any $s,t,u$. It follows that the path
algebra of this subquiver with relations is just $M_n(k)$.
\medskip

Since the arrows of $Q'$ generate $A'$ over the subring
$\times_{i=1}^n k$, the set of paths,
$\{x_{m,i}=e_{1,m}a_{m,i}e_{m+1,1}\}$ for all $m$ and $i$, generate
$A'$ over the subring $M_n(k)$ given by the paths $e_m,f_m$.
They differ from the elements $a_{m,i}$ by multiplication by
invertible paths and we can rewrite our relations between the elements
$a_{m,i}$ and $e_m$ as equivalent relations between the elements
$x_{m,i}$. Moreover they generate the ring $B=e_{1,1}A'e_{1,1}$ and
$A'\cong M_n(B)$. After noting that
$x_{m,i}=e_{1,m}a_{m,i}e_{m+1,n}e_{n,1}$, we see that the relations in
the set $T$ can be rewritten in terms of the elements $x_{m,i}$ as
$x_{m,i}-x_{1,i}$; therefore, we write $x_i=x_{m,i}$ and find the
relations between these elements induced by the relations in $Y'$. We
note that a word $a_{1,i_1}\dots a_{u,i_{u}} = x_{i_1}\dots
x_{i_{u}}e_{1,u+1}$ and so taking a relation $y'_j$ in $Y'$, we see
that the corresponding relation between the elements $x_i$ is
$y_j$. Thus $B$ is isomorphic to $S$ as required.
\end{proof}

\begin{proposition}
\label{p0} The path algebra $A$ in Theorem 1.1 has 
global dimension $\leqslant 2$. 
\end{proposition}
\begin{proof} We consider the homological
dimension of the simple representations of the quiver. There is one
simple $S_i$ for each vertex of the quiver; this is the representation
which assigns the field $k$ to the vertex $i$ and $0$ to every other
vertex and where each arrow gives the zero map. The simple representation
$S_n$ is also the projective representation $P_n$. Because there are no
relations on the full
subquiver on the vertices $\{2,\dots,n\}$, the simple representations
$S_m$ for $m = 2$ to $n-1$ are of homological dimension $1$; in fact,
for $m = 2$ to $n-1$ we have a short exact sequence
\[
\CD
0 @>>> \bigoplus\limits_{i=0}^a P_{m+1} @>{\phi}>> P_m @>>> S_m @>>> 0
\endCD
\]
where the $0$th component of $\phi$ is $\hat{e}_m$ and the $i$th
component is $\hat{a}_{mi}$ for $i>0$.
Now it is clear that the simple $S_1$ has homological dimension $2$
since the kernel of the homomorphism from $P_1$ to $S_1$ has only
simples of the form $S_m$ for $m>1$ as composition factors so that it
has homological dimension at most $1$ and it is not projective since
it is a factor of $\oplus_{i=0}^a P_2$ by a semisimple
subrepresentation.
\end{proof}

Of course, there is no reason to suppose that $A$ is a subalgebra of
$M_n(S)$. Any relation in $S$ between the elements $x_i$ such that
the longest monomial has length less than $n-1$ will give nonzero
elements of $A$ whose image is $0$ in $M_n(S)$. However, the image,
$\bar{A}$ of $A$ in $M_n(S)$ is the path algebra of a quiver with
relations on the same vertex and arrow set so the quiver is directed
and $\bar{A}$ must have finite global dimension. Moreover, it is clear
that $\bar{A}_{\sigma}$ is isomorphic to $\sigma^{-1}A$. In fact, it is
a fairly simple matter to describe $\bar{A}$. We consider the
filtration of $S$ induced by saying the generators have degree $1$;
that is, $S_{0}=k$ and $S_i$ is the finite dimensional vector space
spanned by the monomials in the generators of length at most $i$. Then
$\bar{A}$ is the upper triangular subalgebra of $M_n(S)$ whose
elements have entries from $S_i$ in the $i$th diagonal where the
main diagonal is taken to be the $0$th.
\medskip

We summarise this in the following theorem.

\begin{theorem}
\label{t1}
Let $S$ be a finitely presented algebra. Let the largest degree of a
relation be $n-1$. Then there is an upper triangular finite dimensional
subalgebra $C$ of $M_n(S)$ of which $M_n(S)$ is a universal
localization. In particular, $C$ has finite global dimension.
\hfill\qed
\end{theorem}

In order to see that the examples in the last lemma usually give us
examples of universal localizations that are not stably flat we note
the following lemma.

\begin{lemma}\label{l0}
If $\phi\colon R \to S$ is a stably flat epimorphism of rings then
$$\hbox{\rm global dimension}(S) \leqslant \hbox{\rm global dimension}(R)~.$$
\end{lemma}
\begin{proof}
That $\phi$ is an epimorphism of rings is equivalent to the condition
that the multiplication map from $S\otimes_R S$ to $S$ is an
isomorphism (\cite{Knight70}).  Therefore, by Lemma 3.30 of \cite{Neeman-Ranicki01},
$\Tor_i^R(S,M)=0$ for any $S$ module $M$ and $S\otimes_R M = M$.
Therefore, we can construct a projective resolution of $M$ by applying
$S\otimes_R\_$ to a projective resolution of $M$ as $R$ module.  It
follows that the homological dimension of $M$ as $S$ module is bounded
by its homological dimension as $R$ module.
\end{proof}

There are many possible variations on this method for representing algebras
as universal localizations of finite dimensional algebras. We give two examples
to illustrate possible changes. Let $Q$ be the quiver with relations
having vertices $1,2,3,4$ and arrows $e_1,x_1$ from $1$ to $2$,
$e_2,y_2$ from $2$ to $3$ and $e_{3},x_{3}$ from $3$ to $4$
together with relations $x_1e_2e_{3}-e_1e_2x_{3}$ and
$x_1y_2e_{3}-e_1y_2x_{3}-e_1e_2e_{3}$. On inverting the
arrows $e_1,e_2,e_{3}$ the path algebra we obtain is
$M_{4}(R_1)$ where $R_1$ is the first Weyl algebra.
\medskip

Let $Q$ be the quiver with relations having vertices $1,2,3,4$ and
arrows $e_1,x_1$ from $1$ to $3$, $e_2,y_2$ from $2$ to $3$
and $e_{3},x_{3}$ from $3$ to $4$ together with relations
$e_1x_{3}-x_1e_{3}$, $x_1x_{3}$ and $y_2x_{3}$. On inverting
the arrows $e_1,e_2,e_{3}$ the path algebra we obtain is
$M_{4}(k\langle x,y: x^2, yx\rangle )$. The important point in this
example is that the set of arrows we invert can be simply a maximal
subtree of the quiver, and there may occasionally be an advantage to
doing this if the relations we are interested in can be described
compactly on a tree.
\medskip

At this point, we return to the question of whether there can be an
algorithm to determine the equality of elements in a universal
localization. Thus let $A$ be a ring and $\sigma$ a set of maps
between finitely generated projective modules over $A$. The Malcolmson
normal form states that every element of the localised ring
${\sigma^{-1}A}$ can be written in the form $as^{-1}b$ where $s\colon
P\rightarrow Q$ lies in the upper triangular closure of $\sigma$,
$a\colon A\rightarrow Q$ and $b\colon P\rightarrow A$ are maps in the
category of finitely generated projective modules over the original
ring $A$ and gives an equivalence condition on such elements which
determines when they define the same element of the localised
ring. This equivalence condition depends on the existence of certain
maps in $\sigma$ and the category of finitely generated projective
modules. One might reasonably ask if such an equation could be
constructed algorithmically. In order to show that this is not
possible we do not need to know the exact nature of the equivalence
relation defined by Malcolmson since it is simply important to be
able to demonstrate that there can be no algorithm to determine the
equality of two such elements. We say that the {\it equality problem}
for $(A,\sigma)$ is {\it solvable} if there is an algorithm to
determine the equality of two such elements in ${\sigma^{-1}A}$.

Our proof that the equality problem is not always solvable comes from
the fact that the word problem for groups is not always solvable. Thus
let $G$ be a finitely presented group with generators
$\{x_{i}:1 \leqslant i \leqslant c\}$ and relations 
$\{r_{j}:1 \leqslant j \leqslant d\}$. 
We obtain a finite presentation of its group algebra $kG$ by
taking as generators $\{x_{i},\bar{x}_{i}:1 \leqslant i \leqslant c\}$ 
and as relations $\{x_{i}\bar{x}_{i}-1:1 \leqslant i \leqslant c\}\cup
\{\bar{x}_{i}x_{i}-1:1 \leqslant i\leqslant c\}\cup \{s_{j}-1:1 \leqslant j \leqslant 
d\}$ where $s_{j}$ is obtained from $r_{j}$ by replacing each
occurrence of each $x_{i}^{-1}$ by $\bar{x}_{i}$. Let $A$ be the
finite dimensional algebra we produce by the method considered in and
preceding theorem \ref{t0} and let $\sigma$ be the set of maps between
finitely generated projective modules over $A$ considered there so
that ${\sigma^{-1}A}$ is isomorphic to $M_{n}(kG)$ for a suitable
integer $n$.

\begin{theorem}
\label{eqnotsol}
Let $G$ be a finitely presented group for which the word problem is
not solvable. Let $A$ be the finite dimensional algebra and let
$\sigma$ be the set of maps between finitely generated projective
modules considered in the previous paragraph. Then the equality
problem for $(A,\sigma)$ is not solvable.
\end{theorem}
\begin{proof}
We use the notation developed before theorem \ref{t0}.
The group ring occurs as the endomorphism ring of $P_{1}\otimes
M_{n}(kG)$. The generators of the group and their inverses occur as
elements of the form $x_{m1}e_{1}^{-1}$. Therefore words in these
elements can be written algorithmically in the form $as^{-1}b$ for
suitable maps $a$ and $b$ between finitely generated projective
modules and $s$ in the upper triangular closure of $\sigma$. If there
were an algorithm to determine whether such an element were equal to
the identity map on $P_{1}$ then we would be able to solve the word
problem for the group $G$. Since we cannot solve the word problem,
there can be no such algorithm.
\end{proof}

\section{An explicit computation}
\label{example1}

\ntn{NAi.1}
In this section, let $k$ be a ring and $S$ a $k$-ring, i.e. a ring
homomorphism $k \to S$. We will assume throughout that $S$ is flat as a left $k$-module.
\entn

We define a functor from the category of left $S$-modules to itself.

\dfn{DAi.2}
Recall the short exact sequence of $S$ bimodules
\begin{equation*}
\CD
0 @>>> \Omega_k(S) @>>> S\otimes_k S @>m>> S @>>> 0
\endCD
\end{equation*}
where $\Omega_k(S)$ is the universal bimodule of derivations of $S$ over $k$
and $m$ is the multiplication map. It is split considered as a sequence of left or right $S$
modules since $S$ is a projective module but not as a sequence of bimodules.

Given a left $S$-module $M$, we define $\ck(M)$ from the exact sequence obtained by tensoring with $M$
on the right
\[
\CD
0 @>>> \ck(M) =\Omega_k(S)\otimes_S M @>>> S\otimes_k M @>\mu_M>> M @>>> 0
\endCD
\]
with $\mu_M=m\otimes_S 1_M$ the multiplication map. Thus $\ck(M)$ is the kernel of the multiplication map
and is isomorphic to $\Omega_k(S)\otimes_S M$.

\edfn

\lem{LAi.3}
As in Definition~\ref{DAi.2}, let $M$ be a left $S$-module.
If $M$ is flat as a (left) $k$-module, then so is $\ck(M)$.
\elem

\prf
By hypothesis, both $M$ and $S$ are flat as left $k$-modules.
It follows that $S\otimes_k^{}M$ is also a
flat left $k$-module. In the exact sequence
\[
\CD
0 @>>> \ck(M) @>>> S\otimes_k^{}M @>\mu_M^{}>> M @>>> 0
\endCD
\]
we now know that both $M$ and $S\otimes_k^{}M$ are flat as left
$k$-modules. From the exact sequence for Tor it now follows that
so is $\ck(M)$.
\eprf

Let $M$ be a left $S$-module, flat over $k$. The above produces for us
exact sequences of left $S$-modules, all flat over $k$
\[
\CD
0 @>>> \ck(M) @>i_M^{}>> S\otimes_k^{}M @>\mu_M^{}>> M @>>> 0 \\
0 @>>> \ck^2(M) @>i_{\ck(M)}^{}>> S\otimes_k^{}\ck(M)
@>\mu_{\ck(M)}^{}>> \ck(M) @>>> 0\\
0 @>>> \ck^3(M) @>i_{\ck^2(M)}^{}>>
S\otimes_k^{}\ck^2(M) @>\mu_{\ck^2(M)}^{}>>
\ck^2(M) @>>> 0
\endCD
\]
Splicing these short exact sequences, we deduce

\lem{LAi.4}
Let $M$ be a left $S$-module, flat over $k$. To make the notation
work nicely, define $\ck^0(M)=M$. For $n\geqslant1$ we have defined
$\ck^n(M)$ above.
For each $j \geqslant 1$
there is an exact sequence of left $S$-modules, all flat over $k$
\[
\CD
0 @>>> \ck^j(M) @>i_{\ck^{j-1}(M)}^{}>>
S\otimes_k^{}\ck^{j-1}(M)  @>>> \cdots \\
\cdots @>>> S\otimes_k^{}\ck(M)
@>i_M^{}\mu_{\ck(M)}^{}>> S\otimes_k^{}\ck^0(M)  @>\mu_M^{}>> M @>>> 0.
\endCD
\]
$\,$\hfill{$\Box$}
\elem

The case of most interest to us is where $M=S$. We can assemble the first
$n$ of these exact sequences in vector form.

\lem{LAi.5}
We have an exact sequence
\[
\CD
\left(\begin{array}{c}
0\\ 0 \\ 0\\  \vdots \\ 0 \\ 0 \\ 0
\end{array}\right) \ri\left(\begin{array}{c}
0\\ 0 \\ 0\\  \vdots \\ 0 \\ 0 \\ k
\end{array}\right)\otimes_k^{} \ck^{n-1}(S)  \ri
\left(\begin{array}{c}
0\\ 0 \\ 0\\  \vdots \\ 0 \\ k \\ S
\end{array}\right)\otimes_k^{} \ck^{n-2}(S)  \ri
\cdots
 \endCD
\]
\[
\CD
\cdots @>>>
\left(\begin{array}{c}
k\\ S \\ S\\  \vdots \\ S \\ S \\ S
\end{array}\right)\otimes_k^{} \ck^0(S)  @>>>
\left(\begin{array}{c}
S\\ S \\ S\\  \vdots \\ S \\ S \\ S
\end{array}\right) @>>>
\left(\begin{array}{c}
0\\ 0 \\ 0\\  \vdots \\ 0 \\ 0 \\ 0
\end{array}\right) .
\endCD
\]
$\,$\hfill{\qed}
\elem

\dfn{DAi.6}
Let $A$ be the ring of $n\times n$ lower triangular matrices
\[
A\eq
\left(\begin{array}{ccccccccccc}
k & 0 & 0& \cdots & 0 & 0 & 0\\
S & k & 0& \cdots & 0 & 0 & 0\\
S & S & k& \cdots & 0 & 0 & 0\\
\vdots & \vdots & \vdots&  & \vdots & \vdots & \vdots\\
S & S & S& \cdots & k & 0 & 0\\
S & S & S& \cdots & S & k & 0\\
S & S & S& \cdots & S & S & k
\end{array}\right)
\]
That is, the terms above the diagonal vanish, the diagonal terms lie in $k$,
while the terms below the diagonal may be any elements of $S$.
\edfn

The columns of the matrix ring $A$
are left $A$-modules. We denote them
\[
P_1=\left(\begin{array}{c}
k\\ S \\ S\\  \vdots \\ S \\ S \\ S
\end{array}\right),\quad
P_2=\left(\begin{array}{c}
0\\ k \\ S\\  \vdots \\ S \\ S \\ S
\end{array}\right),\qquad\cdots\qquad
P_{n-1}=\left(\begin{array}{c}
0\\ 0 \\ 0\\  \vdots \\ 0 \\ k \\ S
\end{array}\right),\quad
P_n=\left(\begin{array}{c}
0\\ 0 \\ 0\\  \vdots \\ 0 \\ 0 \\ k
\end{array}\right)
\]
and as a left $A$-module
$$A~=~P_1 \oplus P_2 \oplus \dots \oplus P_n~.$$
Then Lemma~\ref{LAi.5} says that we have an exact sequence
\[
\CD
0 @>>>P_n \otimes_k^{} \ck^{n-1}(S)  @>>>
P_{n-1}\otimes_k^{} \ck^{n-2}(S)  @>>>
\cdots
 \endCD
\]
\[
\CD
\cdots @>>>
P_1\otimes_k^{} \ck^0(S)  @>>>
\left(\begin{array}{c}
S\\ S \\ S\\  \vdots \\ S \\ S \\ S
\end{array}\right) @>>>
\left(\begin{array}{c}
0\\ 0 \\ 0\\  \vdots \\ 0 \\ 0 \\ 0
\end{array}\right) .
\endCD
\]
This is clearly a resolution of left $A$-modules. The modules $P_i$ are all
direct summands of $A$, hence they are projective left $A$-modules. Being
projective, they are certainly flat left $A$-modules. The modules
$\ck^{i-1}(S)$ are flat left $k$-modules. It follows that $P_i\otimes_k^{}
\ck^{i-1}(S)$ are all flat left $A$-modules. Summarizing the
above, we have

\lem{LAi.6.5} The left $A$-module $$N~=~\left(\begin{array}{c}
S\\ S \\ S\\  \vdots \\ S \\ S \\ S
\end{array}\right)$$
has a flat resolution
\[
\CD
0 \ri P_n \otimes_k^{} \ck^{n-1}(S)@>>> \dots  @>>> P_1\otimes_k\ck^0(S)
@>>> N \ri 0~.
\endCD
\]
$\,$\hfill{\qed}
\elem

Define also the right $A$-module
$$M~=~\left(\begin{array}{ccccccc}
S& S & S&  \cdots & S & S & S
\end{array}\right)~.$$

\lem{LAi.7} We have
$$M\otimes_AP_i~=~S~~( 1 \leqslant  i \leqslant  n)$$
and
$$M_n(S)\otimes_AP_i~=~N~~( 1 \leqslant  i \leqslant n)~.$$
\elem
\prf We begin with $M\otimes_AP_i=S$. There are obvious maps
$$\CD
S @>\alpha_i>>M\otimes_AP_i @>\beta_i>> S
\endCD$$
defined by
$$\begin{array}{l}
\alpha_i(s)~=~(0,\dots,0,s)\otimes \begin{pmatrix}
0 \\ \vdots \\ 0 \\ 1 \end{pmatrix}~,\\
\beta_i\bigg((s_1,s_2,\dots,s_n)\otimes\begin{pmatrix}
0 \\ \vdots \\ 0 \\ x_i \\ \vdots \\ x_n \end{pmatrix}\bigg)~=~
\sum\limits^n_{j=i}s_jx_j~.
\end{array}$$
It is clear that the composite $\beta_i\alpha_i$ is the identity.
It suffices to show that $\alpha_i$ is surjective, which we leave to the reader.\\
\indent  The identity $M_n(S)\otimes_AP_i=N$ reduces to the above, after
observing that $M_n(S)=\bigoplus^n_{i=1}M$ as a right $A$-module.
\eprf

\pro{LAi.8} The $\Tor$-groups are
$$\Tor^A_i(M,N)~=~\begin{cases}
S&if~\hbox{$i=0$}\cr
\ck^n(S)&if~\hbox{$i=n-1$}\cr
0&otherwise~.
\end{cases}$$
Consequently
$$\Tor^A_i(M_n(S),M_n(S))=M_n(\Tor^A_i(M,N))=
\begin{cases}
M_n(S)&if~\hbox{$i=0$}\cr
M_n(\ck^n(S))&if~\hbox{$i=n-1$}\cr
0&otherwise~.
\end{cases}$$
\epro
\prf By definition, $\Tor^A_i(M,N)$ is the $i$th homology of the
complex obtained from any flat resolution of $N$ by tensoring over $A$
with $M$. We use the resolution provided by Lemma \ref{LAi.6.5}.
Lemma \ref{LAi.7} allows us to identify $\Tor^A_i(M,N)$ with the $i$th
homology of
$$
S\otimes_k\ck^{n-1}(S) \ri
S\otimes_k^{}\ck^{n-2}(S) \ri \cdots
\ri S\otimes_k^{}\ck^1(S) \ri S\otimes_k^{}\ck^0(S)~.$$
\eprf

\dfn{DAi.9}
Let $\phi:k\ri S$ be the ring homomorphism giving $S$ the structure of
an $A$-ring.
Define $\s$ to be the set of maps $s_i:P_n\ri P_i$ given by the matrices
\[
\CD
\left(\begin{array}{c}
0\\
0\\
\vdots\\
0\\
0\\
\vdots\\
0\\
\phi
\end{array}\right)\quad : \quad\left(\begin{array}{c}
0\\ 0 \\   \vdots \\ 0 \\ 0 \\ \vdots \\ 0 \\ k
\end{array}\right) @>>>
\left(\begin{array}{c}
0\\ 0 \\   \vdots \\ k \\ S \\ \vdots \\ S \\ S
\end{array}\right)
\qquad
\endCD
\]
\edfn

\lem{LAi.10} The ring homomorphism $A \to M_n(S)$ is $\s$-inverting.
\elem
\prf By Lemma \ref{LAi.7}
$$1\otimes s_i~:~M_n(S)\otimes_AP_n \to M_n(S)\otimes_AP_i$$
can be identified with $1:N \to N$.
\eprf

\thm{TAi.11} For $n \geqslant 3$, $A \to M_n(S)$ is universally $\s$-inverting,
$$\s^{-1}A~=~M_n(S)~.$$
\ethm
\prf Let $T$ be a $\s$-inverting $A$-ring. We need to exhibit a
unique factorization
$$A \to M_n(S) \to T~.$$
It follows from $A=\bigoplus\limits^n_{i=1}P_i$ that
$$T~=~\bigoplus\limits^n_{i=1}T\otimes_AP_i$$
with the $T\otimes_AP_i$'s isomorphic f.g. projective $T$-modules. Also,
$$T~=~{\rm End}_T(T)~=~M_n({\rm End}_T(T\otimes_AP_1))~.$$
It therefore suffices to produce a homomorphism
$$S \to {\rm End}_T(T\otimes_AP_1)~.$$
For $x \in S$ define the $A$-module morphisms
\[
\CD
r_x~:~P_i~=~\left(\begin{array}{c}
0\\ \vdots \\ 0 \\ \vdots \\  k \\  \vdots \\ S
\end{array}\right) @>>>
P_j~=~\left(\begin{array}{c}
0\\  \vdots \\ k \\ \vdots \\ S \\ \vdots \\ S
\end{array}\right)
\qquad
\endCD
\]
with components right multiplication by $x$. Define
$$S \ri {\rm End}_T(T\otimes_AP_1)$$
by sending $x \in S$ to
$$\CD P_1 @>(r_1)^{-1}>> P_2 @>r_x>> P_1 \endCD~.$$
Because $r_{x+y}=r_x+r_y$ this is a homomorphism of abelian groups.
The multiplicative identity $r_{xy}=r_xr_y$ follows from the commutative
diagram
$$\CD
P_n @>r_1>> P_{n-1} @>r_1>> P_1\\
@V r_y VV @VV r_y V\\
P_{n-1} @>r_1>> P_1 \\
@V r_x VV \\
P_1
\endCD$$
(This diagram only makes sense if $n > n-1 >1$, i.e. $n \geqslant 3$).
\eprf

\rmk{AAi.12} Suppose $k$ is a field. For any finite-dimensional $k$-algebra
$S$ let $d={\rm dim}_k(S)$. It follows from the exact sequence
$$0 \ri \ck(M) \ri S\otimes_kM \ri M \ri 0$$
that for any f.g. $S$-module $M$
$${\rm dim}_k\ck(M)~=~(d-1){\rm dim}_k(M)~.$$
By induction, for $M=S$ and $n \geqslant 1$
$${\rm dim}_k\ck^n(S)~=~(d-1)^nd~.$$
Thus if $n \geqslant 3$ and $d>1$ then
$$\Tor^A_{n-1}(\sigma^{-1}A,\sigma^{-1}A)~=~M_n(\ck^n(S)) \neq 0~.$$
In particular, for $S=k[\varepsilon]/(\varepsilon^2)$ ($d=2$) the ring
$$A~=~\left(\begin{array}{ccccccccccc}
k & 0 & 0\\
S & k & 0\\
S & S & k
\end{array}\right)$$
is the path algebra of the quiver with relations $(Q,R)$ constructed as
in section 1, with $\sigma^{-1}A=M_3(S)$,
$\Tor^A_2(\sigma^{-1}A,\sigma^{-1}A)=M_3(S)$.  The global dimension of
$S$ is $\infty$, while the global dimension of $A$ is $2$ (cf. 
Proposition 1.2) so $M_3(S)=\sigma^{-1}A$ can also be seen to be non
stably flat by Lemma 1.4.\\ \hskip150pt{\hphantom{AA}}
\ermk

\providecommand{\bysame}{\leavevmode\hbox to3em{\hrulefill}\thinspace}
\providecommand{\MR}{\relax\ifhmode\unskip\space\fi MR }
\providecommand{\MRhref}[2]{%
  \href{http://www.ams.org/mathscinet-getitem?mr=#1}{#2}
}
\providecommand{\href}[2]{#2}

\end{document}